\documentclass{article}

\usepackage[T1]{fontenc}
\usepackage[utf8]{inputenc}

\usepackage{amsthm}
\usepackage{amssymb}
\usepackage{amsmath}
\usepackage{nicefrac}
\usepackage{graphicx}
\usepackage{float}
\usepackage{url}

\usepackage{color}

\title{Numerical Semigroups generated by Primes}
\date{\today}
\author{M. Hellus\and A. Rechenauer\and R. Waldi}
\begin{document}
\maketitle

\begin{abstract}
Let $p_1=2, p_2=3, p_3=5, \ldots$ be the consecutive prime numbers, $S_n$ the numerical semigroup generated by the primes not less than $p_n$ and $u_n$ the largest irredundant generator of $S_n$. We will show, that

\begin{itemize}

\item $u_n\sim3p_n$.

\end{itemize}

Similarly, for the largest integer $f_n$ not contained in $S_n$, by computational evidence (\cite{table1}) we suspect that

\begin{itemize}

\item $f_n$ is an odd number for $n\geq5$ and

\item $f_n\sim3p_n$; further

\item $4p_n>f_{n+1}$ for $n\geq1$.

\end{itemize}

If $f_n$ is odd for large $n$, then $f_n\sim3p_n$. In case $f_n\sim3p_n$ every large even integer $x$ is the sum of two primes. If $4p_n>f_{n+1}$ for $n\geq1$, then the Goldbach conjecture holds true.

Further, Wilf's question in \cite{wilf78} has a positive answer for the semigroups $S_n$.

\end{abstract}
\renewcommand{\thefootnote}{}
\footnote{Michael Hellus, Fakult\"at f\"ur Mathematik, Universit\"at Regensburg, Germany, michael.hellus@mathematik.uni-regensburg.de}
\footnote{Anton Rechenauer, antonrechenauer@gmail.com}
\footnote{Rolf Waldi, Fakult\"at f\"ur Mathematik, Universit\"at Regensburg, Germany, rolf.waldi@mathematik.uni-regensburg.de}



\section*{Introduction}

A \emph{numerical semigroup} is an additively closed subset $S$ of $\mathbb{N}$ with $0\in S$ and only finitely many positive integers outside from $S$, the so-called \emph{gaps} of $S$. The \emph{genus} $g$ of $S$ is the number of its gaps. The set $E=S^*\setminus(S^*+S^*)$, where $S^*=S\setminus\{0\}$, is the (unique) minimal system of generators of $S$. Its elements are called the \emph{atoms} of $S$; their number $e$ is the \emph{embedding dimension} of $S$. The \emph{multiplicity} of $S$ is the smallest element $p$ of $S^*$.

From now on we assume that $S\neq\mathbb{N}$. Then the greatest gap $f$ is the \emph{Frobenius number} of $S$. Since $(f+1)+\mathbb{N}\subseteq S^*$ we have $(p+f+1)+\mathbb{N}\subseteq p+S^*$, hence the atoms of $S$ are contained in the interval $[p,p+f]$.

For our investigation of certain numerical semigroups $S$ generated by prime numbers, the fractions
\[\frac f p\text{, }\frac{1+f}p\text{, }\frac g {1+f}\text{ and }\frac{e-1}e\]
will play a role. For general $S$, what is known about these fractions?

First of all it is well known and easily seen that
\[\frac12\leq\frac g{1+f}\leq\frac{p-1}p,\]
and both bounds for $\frac g{1+f}$ are attained.

However, the following is still open:

\vspace{.2cm}

\noindent\textbf{Wilf's question} (\cite{wilf78}): Is it (even) true that
\[\tag{1}\frac g{1+f}\leq\frac{e-1}e\]
for every numerical semigroup?

\vspace{.2cm}

A partial answer is given by the following result of Eliahou:

\vspace{.2cm}

\noindent\textbf{\cite[Corollary 6.5]{eliahou18}} If $\frac{1+f}p\leq3$, then $\frac g{1+f}\leq\frac{e-1}e$.

\vspace{.2cm}

In \cite{zhai13}, Zhai has shown that $\frac{1+f}p\leq3$ holds for almost all numerical semigroups of genus $g$ (as $g$ goes to infinity).

Therefore, for randomly chosen $S$, one has $\frac g{1+f}\leq\frac{e-1}e$ almost surely.

\vspace{.2cm}

We shall consider the following semigroups: Let $p_1=2$, $p_2=3$, $p_3=5, \ldots $ be the sequence of prime numbers in natural order and let $S_n$, for $n\geq1$, be the numerical semigroup generated by all prime numbers not less than $p_n$; the multiplicity of $S_n$ is $p_n$ and we denote the aforementioned invariants of $S_n$ by $g_n$, $f_n$, $e_n$ and $E_n$. Since $S_{n+1}$ is a subsemigroup of $S_n$ it is clear that $f_n\leq f_{n+1}$ for all $n\geq1$. The atoms of $S_n$ are contained in the interval $[p_n,p_n+f_n]$; conversely, each odd integer from $S_n\cap[p_n,3p_n[$ is an atom of $S_n$.

As a major result we will see that Wilf's question has a positive answer for $S_n$. Further ${g_n}/{p_n}$ converges to $5/2$ for $n\to\infty$.

The prime number theorem suggests that there should be -- like for the sequence $(p_n)$ -- some asymptotic behavior of $(g_n)$, $(f_n)$ and $(e_n)$.

Based on the list $f_1, f_2, \ldots, f_{2000}$ from \cite{oeisA180306}, extensive calculations (cf. our table 1 in \cite{table1}) gave evidence for the following three conjectures:

\vspace{.2cm}

\begin{enumerate}

\item[(C1)] $f_n\sim 3p_n$, i.\,e. $\lim_{n\to\infty}\frac{f_n}{p_n}=3$,

\end{enumerate}

as already observed by Kl\o{}ve \cite{klove74}, see also the comments in \cite[p.\,56]{erdos_graham_80}; note that Kl\o{}ve works with \emph{distinct} primes, therefore his conjecture is formally stronger than ours, however see also \cite[comment by user ``Emil Je\u{r}\'{a}bek'', Apr 4 '12]{MO93002}.

\vspace{.2cm}

\noindent By Proposition 1, we know that
\[\tag{2}3p_n-f_n\leq6.\]

\begin{enumerate}

\item[(C2)] $f_{n+1}<4p_n$ for all $n\geq1$.

\end{enumerate}

and

\begin{enumerate}

\item[] $3p_n<f_{n+1}$ for $n\geq3$.

\end{enumerate}

\noindent It is immediate from (2) that at least

\[3p_n\leq f_{n+1}\text{ for }n\geq2.\]

\begin{figure}[H]
\includegraphics[width=250px,height=200px,trim=0 40 0 20,clip]{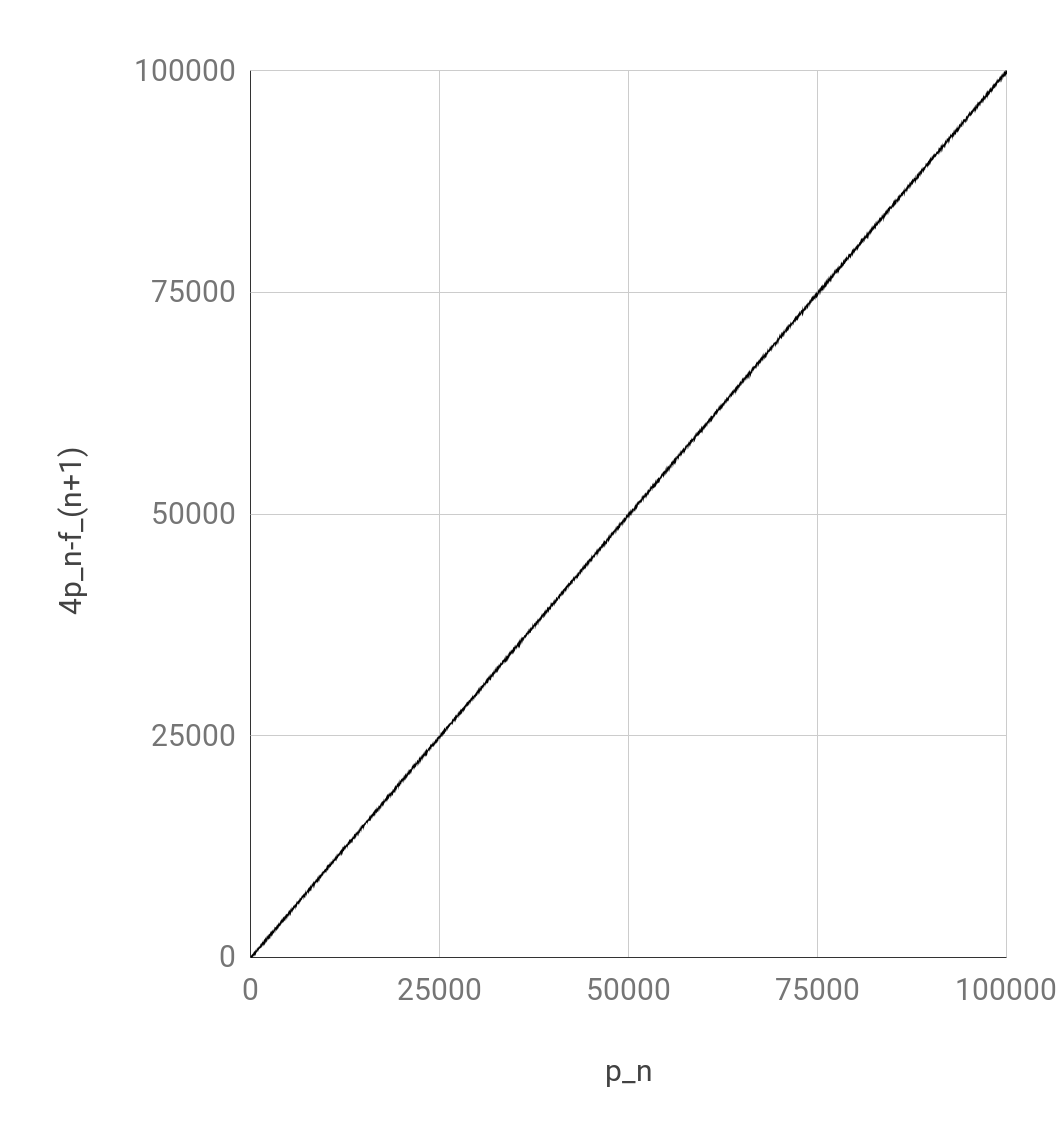}
\caption{$4p_n-f_{n+1}$ vs $p_n$}
\end{figure}

\begin{figure}[H]
\includegraphics[width=250px,height=200px,trim=0 40 0 20,clip]{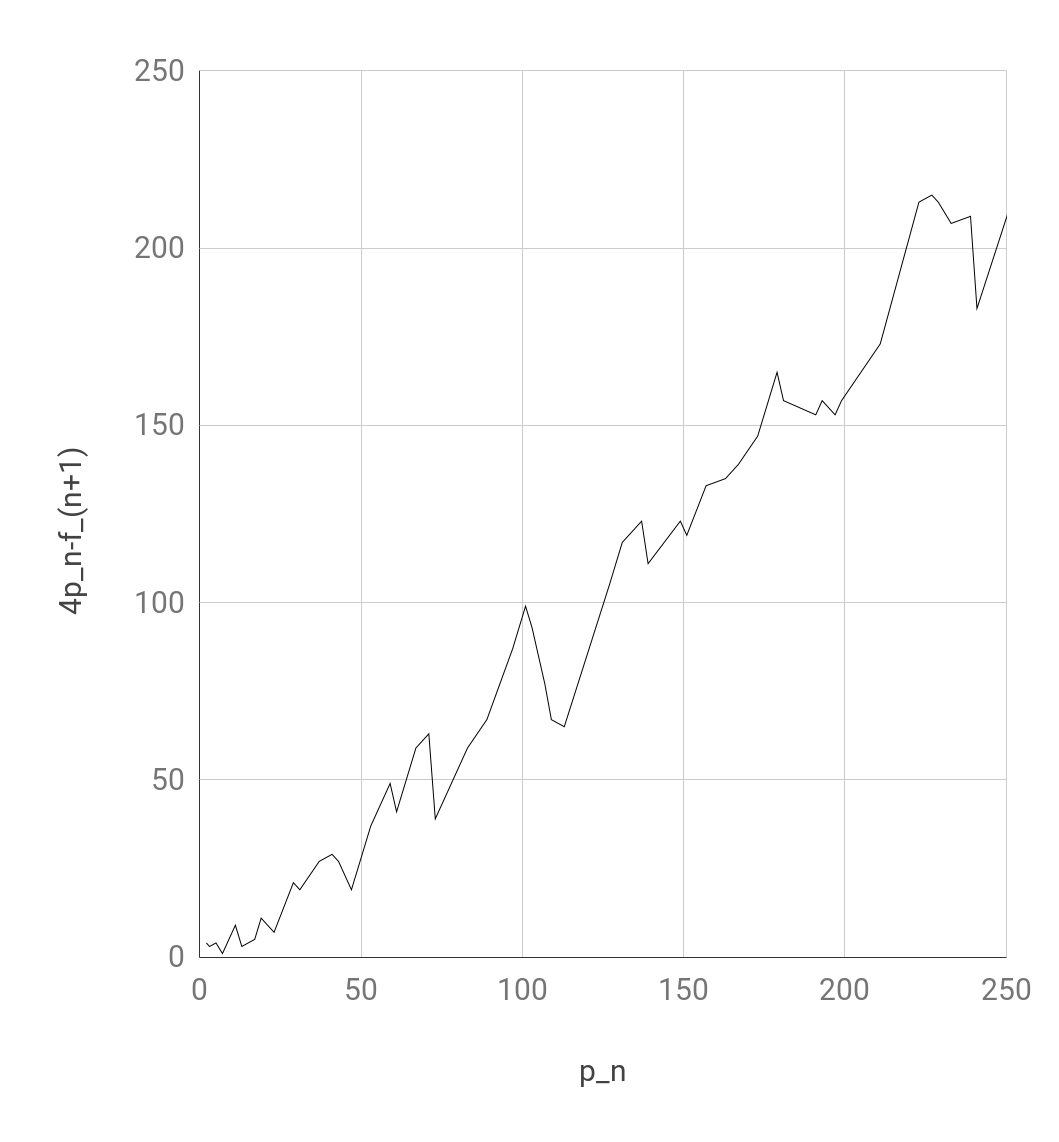}
\caption{$4p_n-f_{n+1}$ vs $p_n$}
\end{figure}

\vspace{.2cm}

As already noticed in \cite{klove74} and in \cite[answer by user ``Woett'', Apr 3 '12]{MO93002}, both conjectures (C1) and (C2) are closely related to Goldbach's conjecture. As we will see in Proposition 4, (C1) is a consequence of conjecture

\begin{enumerate}

\item[(C3)] $f_n$ is odd for $n\geq5$.

\end{enumerate}

Notice again, that a conjecture similar to (C3) was already formulated in \cite{klove74}, however for the (related) notion 'threshold of completeness' for the sequence of all prime numbers, in the sense of \cite{erdos_graham_80}.

\vspace{.2cm}

\noindent Figure 1 indicates, that $\lim_{n\to\infty}\frac{f_n}{p_n}=3$ should be true.

\vspace{.2cm}

\noindent As for (C2), by figure 1 and figure 2, evidently $4p_n-f_{n+1}$ should stay positive for all time.

\vspace{.2cm}

\noindent\textbf{Observations} Numerical experiments suggest that similiar conjectures can be made if one restricts the generating sequence to prime numbers in a fixed arithmetic progression $a+kd$ for $(a,d)=1$. In such a case the limit of $\frac{f_n}{p_n}$ would apparently be $d+1$ ($d$ even)  or $2d+1$ ($d$ odd), see figure 3, and table 2 in \cite{table2}.

\begin{figure}[H]
\centering
\includegraphics[width=300px,trim=0 0 0 0,clip]{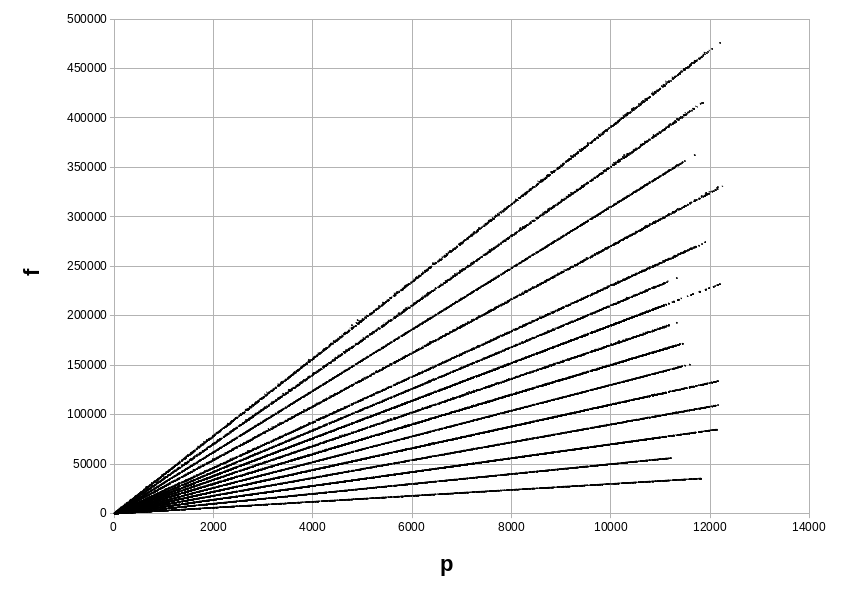}
\caption{$f$ vs. $p$ for some series of semigroups as in the 'Observations'}
\end{figure}

\vspace{.2cm}

The following version of Vinogradov's theorem is due to Matom\"aki, Maynard and Shao. It is fundamental for the considerations in this paper.

\vspace{.2cm}

\noindent\textbf{\cite[Theorem 1.1]{matomaki_etal_17}} Let $\theta>\frac{11}{20}$. Every sufficiently large odd integer $n$ can be written as the sum $n=q_1+q_2+q_3$ of three primes with the restriction
\[\left|q_i-\frac n3\right|\leq n^\theta\text{ for }i=1, 2, 3.\]
Of course we could have used just as well one of the predecessors of this theorem, see the references in \cite{matomaki_etal_17}.

\section{Variants of Goldbach's conjecture}

For $x,y\in\mathbb{Q}$, $x\leq y$ we denote by $[x,y]$ the 'integral interval'
\[[x,y]:=\{n\in\mathbb{Z}|x\leq n\leq y\},\]
accordingly we define $[x,y[$, $]x,y]$, $]x,y[$, $[x,\infty[$.

\vspace{.2cm}

For $x\geq2$ we define $S_n^x$ to be the numerical semigroup generated by the primes in the interval $I_n^x:=[p_n,x\cdot p_n[$ and $f_n^x$ its Frobenius number.

\vspace{.2cm}

A minor step towards a proof of conjecture (C1) is

\vspace{.2cm}

\noindent\textbf{Proposition 1}
\[f_n\geq3p_n-6.\]
In particular for the null sequence $r(n):=6/p_n$ we have
\[\frac{f_n}{p_n}\geq3-r(n)\text{ for every }n\geq1.\]

\vspace{.2cm}

\noindent\textbf{Proof} For $n\geq3$, obviously, the odd number $3p_n-6$ is neither a prime nor the sum of primes greater than or equal to $p_n$, hence $3p_n-6$ is not contained in $S_n$.

\hfill$\square$

\vspace{.2cm}

\noindent\textbf{Remark} A final (major) step on the way to (C1) would be to find a null sequence $l(n)$ such that
\[3+l(n)\geq\frac{f_n}{p_n}.\]

\vspace{.2cm}

\noindent\textbf{Proposition 2} If (C1) is true then every sufficiently large even number $x$ can be written as the sum $x=p+q$ of prime numbers $p, q$.

\noindent\textbf{Addendum} The prime number $p$ can be chosen from the interval $]\frac x4,\frac x2]$.

\vspace{.2cm}

\noindent\textbf{Proof} By the prime number theorem, we have $p_{n+1}\sim p_n$. (C1) implies
\[f_{n+1}\sim3p_{n+1}\sim3p_n,\]
i.\,e.
\[\lim_{n\to\infty}\frac{f_{n+1}}{p_n}=3.\]
In particular, there exists $n_0\ge1$ such that $\frac{f_{n+1}}{p_n}<4$ for all $n\geq n_0$.

It remains to show:

\vspace{.2cm}

\noindent\textbf{Lemma 1} If $n_0\geq1$ is such that $\frac{f_{n+1}}{p_n}<4$ for all $n\geq n_0$ then every even number $x>2$ with $x>f_{n_0}$ can be written as the sum
\begin{equation}\tag{1}x=p+q\text{ with prime numbers }p\leq q\text{ and such that }\frac x4<p\leq \frac x2.\end{equation}

\noindent\textbf{Proof} By our hypothesis,
\[f_n\leq f_{n+1}<4p_n<4p_{n+1}\text{ for all }n\geq n_0\]
and hence, for $I_n:=[1+f_n,4p_n[$ ($n\geq n_0$),
\[[1+f_{n_0},\infty[=\bigcup_{n\geq n_0}I_n.\]
Therefore it suffices to prove (1) for all even numbers $x>2$ from the interval $I_n$, for $n\geq n_0$.

By definition of $f_n$, every $x\in I_n$ can be written as the sum of primes $p\geq p_n$.

If in addition $x>2$ is even, then, because of $f_n<x<4p_n$, the number $x$ is the sum of precisely two prime numbers $p\leq q$ with
\[p_n\leq p\leq q=x-p<4p_n-p\leq3p,\]
hence
\[\frac x4<p\leq\frac x2.\]
\hfill$\square$

The special case $n_0=1$ of Lemma 1 gives

\noindent\textbf{Proposition 3} If (C2) is true then every even number $x>2$ can be written as the sum $x=p+q$ of prime numbers $p\leq q$ as described in the Addendum above. In particular for each $n\geq1$, $4p_n=p+q$ with primes $p_{n+1}\leq p\leq q$, implying Bertrand's postulate.\hfill$\square$

\vspace{.2cm}

\vspace{.2cm}

\noindent\textbf{Proposition 4} If the Frobenius number $f_n$ is odd for all large $n$, then $f_n\sim3p_n$. In particular, conjecture (C3) implies conjecture (C1).

\vspace{.2cm}

\noindent\textbf{Proof}
From \cite[Theorem 1.1]{matomaki_etal_17} we get:

\vspace{.2cm}

\noindent\textbf{Lemma 2} Let $\varepsilon>0$. For odd $N$ large enough, there are prime numbers $q_1$, $q_2$, $q_3$ with
\[N=q_1+q_2+q_3\]
and such that
\[\frac 1{3+\varepsilon}\cdot N<q_i<\frac{3+2\varepsilon}{9+3\varepsilon}\cdot N\text{, i.\,e. }\left|q_i-\frac N3\right|<\frac\varepsilon{9+3\varepsilon}\cdot N\text{ for }i=1,2,3.\]
\noindent\textbf{Proof of Lemma 2} The claim follows immediately from \cite[Theorem 1.1]{matomaki_etal_17}, since $\theta:=\frac35>\frac{11}{20}$ and, for large $N$, $N^\frac35<\frac\varepsilon{9+3\varepsilon}\cdot N$.\hfill$\square_{\text{Lemma 2}}$

\vspace{.2cm}

By our hypothesis, $f_{n+1}$ is odd for large $n$. In Lemma 3 below we will show that, for each $\varepsilon>0$, we have $f_{n+1}<(3+\varepsilon)p_n$ for large $n$; then the claim of Proposition 4 follows from Proposition 1.\hfill$\square_{\text{Proposition 4}}$

\vspace{.2cm}

\noindent\textbf{Lemma 3} Let $\varepsilon>0$. Then for large $n$, each odd integer $N\geq(3+\varepsilon)p_n$ is contained in $S_{n+1}$. In particular, for large $n$
\[f_{n+1}<(3+\varepsilon)p_n\text{ if }f_{n+1}\text{ is odd, and}\]
\[f_{n+1}<(3+\varepsilon)p_n+p_{n+1}\text{ if }f_{n+1}\text{ is even,}\]
since then $f_{n+1}-p_{n+1}$ is odd and not in $S_{n+1}$.

\noindent\textbf{Proof} Since $N$ is odd and large for large $n$, by Lemma 2 there exist prime numbers $q_1$, $q_2$, $q_3$ with
\[N=q_1+q_2+q_3\]
and such that
\[\frac N{3+\varepsilon}<q_i\text{ for }i=1,2,3.\]
By assumption, $\frac N{3+\varepsilon}\geq p_n$, hence
\[q_i>p_n\text{, i.\,e. }q_i\geq p_{n+1}\]
for the prime numbers $q_i$. This implies $N=q_1+q_2+q_3\in S_{n+1}.$\hfill$\square$

\vspace{.2cm}

For a similar argument, see \cite[answer by user ``Anonymous'', Apr 5'12]{MO93002}.

\vspace{.2cm}

\noindent\textbf{Remarks}

\noindent a) It is immediate from Lemma 3 that
\[\limsup_{n\to\infty}\frac{f_n}{p_n}\leq4.\]
As a consequence, a proof of $\limsup_{n\to\infty}\frac{f_n}{p_n}\neq4$ would imply the binary Goldbach conjecture for large $x$ with the Addendum from above -- see Lemma 1 and the proof of Proposition 2.

\vspace{.2cm}

\noindent b) The estimate $\limsup_{n\to\infty}\frac{f_n}{p_n}\leq4$ together with a sketch of proof was already formulated in \cite[comment by user ``Fran\c{c}ois Brunault'' (Apr 6 '12) to answer by user ``Anonymous'' (Apr 5 '12)]{MO93002}. Our proof is essentially an elaboration of this sketch.

\vspace{.2cm}

\noindent c) Lemma 3 shows that
\[f_{n+1}<5p_{n+1}\text{ for large }n.\]
Because of $p_{n+1}<2p_n$ (Bertrand's postulate) this implies also that there exists a constant $C$ with
\begin{equation}\tag{2}f_{n+1}<Cp_n\text{ for all }n.\end{equation}

\noindent Conjecture (C2) says that in (2) one can actually take $C=4$.
\vspace{.2cm}

\noindent Notice that (2) already follows from \cite[Lemma 1]{benkoski_erdos_74}.

\vspace{.2cm}

\noindent\textbf{Problem} Find an explicit pair $(n_0,C_0)$ of numbers such that
\[f_{n+1}<C_0\cdot p_n\text{ for every }n\geq n_0.\]

\vspace{.2cm}

Next we shall study the asymptotic behavior of the set of atoms of $S_n$.

\noindent Lemma 2 will imply

\vspace{.2cm}

\noindent\textbf{Corollary} Let $\varepsilon>0$. Then $S_n=S_n^{3+\varepsilon}$ for large $n$.

\vspace{.2cm}

In particular, $E_n\subseteq[p_n,(3+\varepsilon)p_n[$ for large $n$, and $\log u_n\sim\log p_n$.

\noindent On the other hand, the primes in $[p_n,3p_n[$ are atoms of $S_n$. hence for large $n$, $\pi(3p_n)\leq\pi(u_n)\leq\pi((3+\varepsilon)p_n)$. The prime number theorem yields
\[3n\leq \pi(u_n)\leq(3+\varepsilon)n\text{ for large n}.\]
Consequently we have the following

\vspace{.2cm}

\noindent\textbf{Theorem} $\pi(u_n)\sim3n$, $e_n\sim2n$ and $u_n\sim3p_n$.

\vspace{.2cm}
\noindent\textbf{Proof of the Corollary} It suffices to prove the claim for arbitrarily small values of $\varepsilon$:

First we show that, if $\varepsilon<3$, then
\[S_{n+1}^{3+\varepsilon}\subseteq S_n^{3+\varepsilon}\]
for large $n$. For this it suffices to show that every prime number $p$ on the interval $[p_{n+1},(3+\varepsilon)p_{n+1}[$ is in $S_n^{3+\varepsilon}$:

Firstly, $p\geq p_{n+1}>p_n$.

Now we distinguish two cases:

\begin{enumerate}

\item[I] $p<(3+\varepsilon)p_n$: Then $p\in I_n^{3+\varepsilon}$, hence $p\in S_n^{3+\varepsilon}$.

\item[II] $p\geq(3+\varepsilon)p_n$: For $n$ large enough, by Lemma 2 there exist prime numbers $q_1,q_2,q_3$ with
\[p=q_1+q_2+q_3\]
and such that
\[p_n\buildrel\text{II}\over\leq\frac p{3+\varepsilon}<q_i<\frac{3+2\varepsilon}{9+3\varepsilon}p\text{ for }i=1,2,3.\]
By Chebyshev, Bertrand's postulate $p_{n+1}<2p_n$ holds. Therefore,
\[p\buildrel\text{hypothesis}\over<(3+\varepsilon)p_{n+1}<(6+2\varepsilon)p_n\]
and hence
\[q_i<\frac{3+2\varepsilon}{9+3\varepsilon}p<\frac{3+2\varepsilon}{9+3\varepsilon}(6+2\varepsilon)p_n<(3+\varepsilon)p_n,\]
if $\varepsilon<3$. It follows that
\[q_i\in[p_n,(3+\varepsilon)p_n[\text{ for }i=1,2,3\text{ and hence}\]
\[p=q_1+q_2+q_3\in S_n^{3+\varepsilon},\]
which proves the above claim.

Recursively, we get from $S_{n+1}^{3+\varepsilon}\subseteq S_n^{3+\varepsilon}$ that
\[p_k\in S_k^{3+\varepsilon}\subseteq S_n^{3+\varepsilon}\text{  for all }k\geq n.\]
Therefore,
\[S_n=S_n^{3+\varepsilon}.\]

\end{enumerate}

\hfill$\square$

By \cite[Cor.~6.5]{eliahou18}, for arbitrary numerical semigroups $S$, Wilf's inequality $\frac g{1+f}\leq\frac{e-1}e$ holds, whenever $f<3\cdot p$. Further by \cite{zhai13}, the latter is true for almost every numerical semigroup of genus $g$ (as $g$ goes to infinity).

In contrast, according to table 1 in \cite{table1}, for the semigroups $S_n$, the relation $f_n<3\cdot p_n$ seems to occur extremely seldom, but over and over again (see figure 4).

\begin{figure}[H]
\includegraphics[width=300px,trim=0 0 0 0,clip]{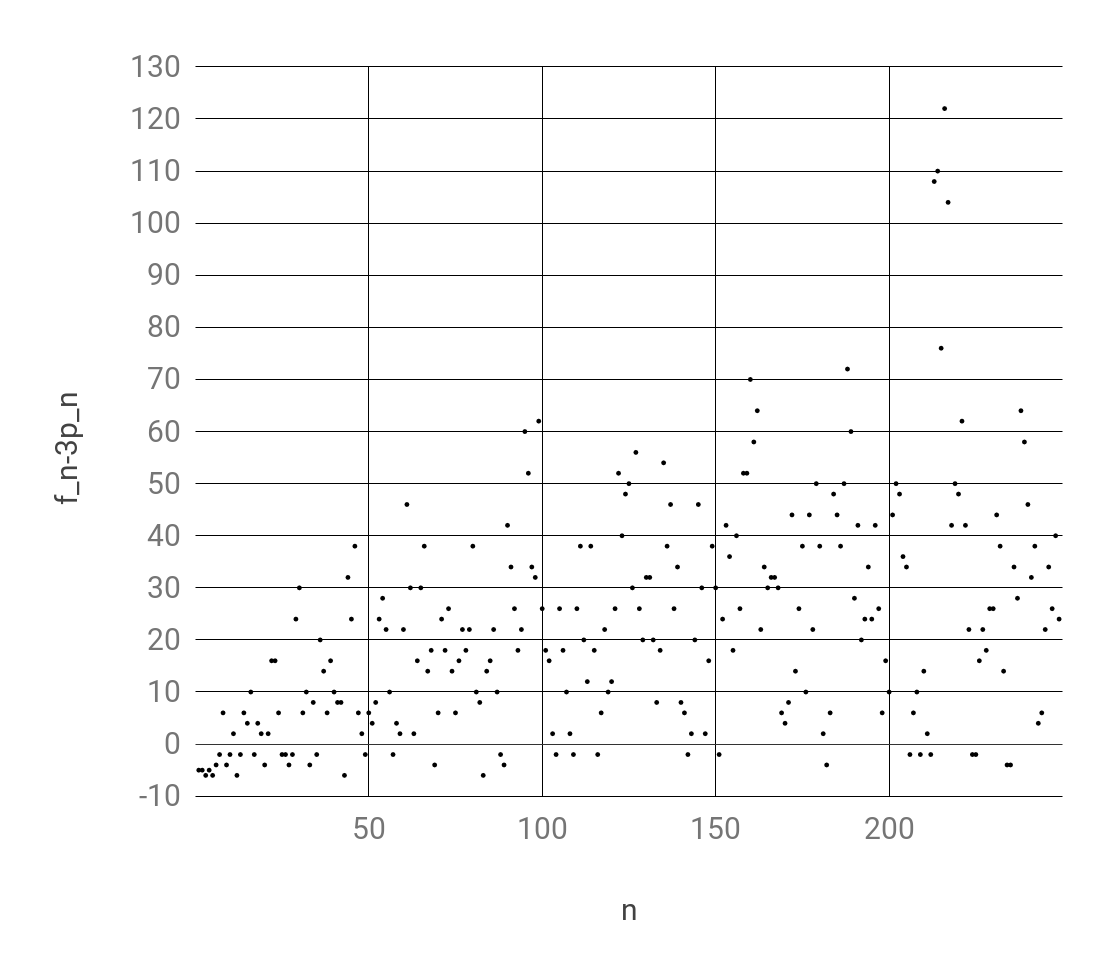}
\caption{$f_n-3p_n$ vs $n$}
\end{figure}

The following considerations are related to \cite[answer by user ``Aaron Meyerowitz'', Apr 3 '12]{MO93002}:

Let $f_n<3\cdot p_n$. Then the odd number $3\cdot p_n+6$ is in $S_n$, but not a prime; hence $p_{n+1}\leq p_n+6$.

\begin{enumerate}

\item[1.] If $p_{n+1}=p_n+4$, since $3\cdot p_n+6\in S_n$ is not a prime, $p_n+6$ must be prime.

\item[2.] If $p_{n+1}=p_n+6$, then the odd numbers $3p_n+2$ and $3p_n+4$ must be atoms in $S_n$, hence primes.

\end{enumerate}

\noindent In any case:

\noindent\textbf{Nota bene} If $f_n<3p_n$, then there is a twin prime pair within $[p_n,3p_n+4]$.

So we cannot expect to prove, that $f_n<3p_n$ happens infinitely often, since this would prove the \textit{twin prime conjecture}, that there are infinitely many twin prime pairs. Another consequence would be that
\[\liminf_{n\to\infty}\frac{f_n}{p_n}=3,\]
since one always has that this limit inferior is $\geq3$, by Proposition 1.

The next section is attended to Wilf's question mentioned above.

\section{The question of Wilf for the semigroups $\mathbf{S_n}$}

\noindent\textbf{Proposition 5} For the semigroups $S_n$, Wilf's (proposed) inequality
\begin{equation}\tag{1}\frac{g_n}{1+f_n}\leq\frac{e_n-1}{e_n}\end{equation}
holds.

\vspace{.2cm}

\noindent\textbf{Proof} For $n<429$, have a look at table 1 in \cite{table1}. Now let $n\geq429$.

\noindent Instead of (1), we would rather prove the equivalent relation
\[\tag{2}e_n(1+f_n-g_n)\geq1+f_n.\]
According to \cite[Cor.~6.5]{eliahou18} we may assume, that $3p_n<1+f_n$. Hence the primes in the interval $[p_n,3p_n[$ are elements of $S_n$ lying below $1+f_n$, and in fact, they are atoms of $S_n$ as well. This implies for the prime-counting function $\pi$
\[\tag{3}e_n(1+f_n-g_n)\geq(\pi(3p_n)-n+1)^2.\]

By Rosser and Schoenfeld \cite[Theorem 2]{rosser_schoenfeld62} we have
\[\tag{4}\pi(x)<\frac x{\log x-\frac32}\text{\ \ for }x>e^{\frac32}\text{, and}\]
\[\tag{5}\pi(x)>\frac x{\log x-\frac12}\text{\ \ for }x\geq67.\]
Further $\lambda(x):=3\cdot\frac{\log x-\frac32}{\log(3x)-\frac12}$ is strictly increasing for $x>1$, hence
\[\tag{6}2n<\pi(3p_n)<3n\text{\ \ for }n\geq429.\]
\noindent\textbf{Proof} Since $\lambda(x)$ is strictly increasing, we get for $n\geq429$, i.\,e. $p_n\geq2971$
\[\pi(3p_n)\buildrel\text{(5)}\over>\frac{3p_n}{\log(3p_n)-\frac12}\buildrel\text{(4)}\over>\pi(p_n)\cdot\lambda(p_n)\geq n\cdot\lambda(2971)>2n\text{, and}\]
\[\pi(3p_n)\buildrel\text{(4)}\over<\frac{3p_n}{\log p_n+\log3-\frac32}<\frac{3p_n}{\log p_n-\frac12}\buildrel\text{(5)}\over<3n\]
\hfill$\square$

In particular, by (3) and (6)
\[e_n(1+f_n-g_n)\buildrel\text{(3)}\over\geq(\pi(3p_n)-n+1)^2\buildrel\text{(6)}\over\geq(n+2)^2.\]
It remains to prove

\vspace{.2cm}

\noindent\textbf{Lemma 4} If $n\geq429$, then
\[f_n<n^2.\]

\vspace{.2cm}

\noindent\textbf{Proof} Let $N\leq a_1<\ldots<a_N$ be positive integers with $(a_1,\ldots,a_N)=1$, $S=\langle a_1,\ldots, a_N\rangle$ the numerical semigroup generated by these numbers and $f$ its Frobenius number. Then, by Selmer \cite{selmer77} we have the following theorem (see the book \cite{ramirezalfonsin_05} of Ram\'{i}rez Alfons\'{i}n). It is an improvement of a former result \cite[Theorem 1]{erdos_graham_72} of Erd\H os and Graham.

\noindent\textbf{\cite[Theorem 3.1.11]{ramirezalfonsin_05}}
\[\tag{7}f\leq 2\cdot a_N\left\lfloor\frac{a_1}N\right\rfloor-a_1.\]
We will apply this to the semigroup $S_n^3\subseteq S_n$ generated by the primes
\[p_n=a_1<p_{n+1}=a_2<\ldots<p_{N+n-1}=a_N\]
in the interval $I_n^3=[p_n,3p_n[$, with Frobenius number $f_n^3$, hence
\[N=\pi(3p_n)-n+1, a_N=p_{\pi(3p_n)}=\text{the largest prime in }I_n^3.\]
By \cite[Theorem 3, Corollary, (3.12)]{rosser_schoenfeld62} we have \[p_n>n\log n\geq n\log429>6n\buildrel\text{(6)}\over>N,\]
hence the above theorem can be applied.

By (6) and (7), $p_{\pi(3p_n)}\buildrel\text{(6)}\over<p_{3n}$ and
\[f_n\leq f_n^3\buildrel\text{(7)}\over<2\cdot p_{\pi(3p_n)}\cdot\frac{p_n}{\pi(3p_n)-n+1}\buildrel\text{(6)}\over<2\cdot p_{3n}\cdot\frac{p_n}{n+2}.\]
From Rosser and Schoenfeld's result \cite[Theorem 3, Corollary, (3.13)]{rosser_schoenfeld62}
\[\tag{8}p_k<k(\log k+\log\log k)\text{ for }k\geq6\]
finally we shall conclude that $2\cdot p_{3n}\cdot\frac{p_n}{n+2}<n^2\text{ for }n\geq429:$

Elementary calculus yields
\[\tag{9}\lambda_2(x):=6\cdot(\log(3x)+\log\log(3x))\cdot(\log x+\log\log x)<x+2\text{ for }x\geq429,\]
since
\[\lambda_2(429)<431\text{ and }\lambda_2'(x)<1=(x+2)'\text{ for }x\geq429.\]
Hence
\[2\cdot p_{3n}\cdot p_n\buildrel\text{(8)}\over<n^2\cdot\lambda_2(n)\buildrel\text{(9)}\over<n^2\cdot(n+2)\text{ for }n\geq429.\]
\hfill$\square$

See also P.~Dusart's th\`{e}se \cite
 {dusart98} for more estimates like (4), (5) and (8).
 
\vspace{.2cm}

\noindent\textbf{Remark} Looking at table 3 in \cite{table3} we see, that even
\[\pi(3p_n)>2n\text{ for }n > 8\text{ and }\pi(3p_n)<3n\text{ for n > 1}\]
(which may be found elsewhere), and
\[f_n\leq n^2\text{ for }n\neq5.\]

\vspace{.2cm}

At last we will see that, apparently, the quotient $\frac{g_n}{1+f_n}$ should converge to $\frac56$ (whereas $\lim_{n\to\infty}\frac{e_n-1}{e_n}=1$, since $e_n\sim2n$ by our Theorem).

\vspace{.2cm}

\noindent\textbf{Proposition 6} The quotient $\frac{g_n}{p_n}$ converges and $\lim_{n\to\infty}\frac{g_n}{p_n}=\frac52$. Hence under the assumption $\lim_{n\to\infty}\frac{p_n}{f_n}=\frac13$ (C1) (which should be true by computational evidence) we have
\[\lim_{n\to\infty}\frac{g_n}{1+f_n}=\frac56.\]
\noindent\textbf{Proof} For that, we consider the proportion $\alpha_k(n)$ of gaps of $S_n$ among the integers in $[k\cdot p_n,(k+1)\cdot p_n]$, ($k,n\geq1$). Besides \cite[Theorem 1.1]{matomaki_etal_17}, we shall need the following similar result about the representation of \emph{even} numbers as the sum of two primes:

\vspace{.2cm}

\noindent\textbf{\cite[Theorem 1, Corollary]{coppola_laporta_95}} Let $\varepsilon>0$ and $A>0$ be real constants. For $N>0$ let $E(N)$ be the set of even numbers $2m\in[N,2N]$, which cannot be written as the sum $2m = q_1 + q_2$ of primes $q_1$ and $q_2$ with the restriction
\[|q_j-m|\leq m^{\frac58+\varepsilon}\text{ for }j=1,2.\]
Then there is a constant $D>0$ such that $\#E(N)<D\cdot N/(\log N)^A$.\hfill$\square$

\vspace{.2cm}

From these two facts together with the prime number theorem, we conclude the following asymptotic behavior of the numbers $\alpha_k(n)$, as $n$ goes to infinity:
\[\alpha_0(n)\to1, \alpha_1(n)\to1, \alpha_2(n)\to\frac12\text{ and }\alpha_k(n)\to0\text{ for }k\geq3.\]
Hence
\[\lim_{n\to\infty}\frac{g_n}{p_n}=1+1+\frac12=\frac52.\]
(Notice that for large $n$, by Lemma 3 we have $f_n<5p_n$, hence $\alpha_k(n)=0$ for $k\geq5$.)\hfill$\square$

\vspace{.2cm}

\noindent\textbf{Remark} Let $f_{n,e}$ be the largest even gap of $S_n$. Our computations (see table 1 in \cite{table1}) suggest that $f_{n,e}\sim 2p_n$. In this case, by Proposition 1 and Proposition 4, $f_n$ is odd for large $n$ and conjecture (C1) holds.

\vspace{.2cm}

\noindent\textbf{Acknowledgement} We thank F.~Brunault, O.~Forster and K.~Matom\"aki for valuable hints.

\end{document}